\documentclass[a4paper,final,intlimits,reqno,10pt,oneside]{amsart}
\usepackage{amssymb,enumerate}

\usepackage[cm,headings]{fullpage}

\newtheorem{theorem}{Theorem}
\newtheorem{lemma}[theorem]{Lemma}

\newcommand\abs[1]{\left|#1\right|}
\newcommand\absb[2]{\csname#1l\endcsname|#2\csname#1r\endcsname|}

\newcommand\norm[1]{\left\|#1\right\|}
\newcommand\normb[2]{\csname#1l\endcsname\|#2\csname#1r\endcsname\|}
\newcommand\norms[1]{\mathopen\|#1\mathclose\|}

\newcommand\de{{\delta}}
\newcommand\ve{{\varepsilon}}
\newcommand\vp{{\varphi}}
\newcommand\al{\alpha}

\newcommand\N{\mathbb N}
\newcommand\R{\mathbb R}
\newcommand\Q{\mathbb Q}
\newcommand\mc{\mathcal}

\renewcommand\d{\,\mathrm{d}}

\begin{document}
\title{On Peano's theorem in Banach spaces}
\author{Petr H\'ajek}
\address{Mathematical Institute\\Czech Academy of Science\\\v Zitn\'a 25\\115 67 Praha 1\\Czech Republic}
\email{hajek@math.cas.cz}
\author{Michal Johanis}
\address{Department of Mathematical Analysis\\Charles University\\Sokolovsk\'a 83\\186 75 Praha 8\\Czech Republic}
\email{johanis@karlin.mff.cuni.cz}
\date{June 2009}
\thanks{Supported by grants: the research project MSM 0021620839, Inst. Res. Plan AV0Z10190503, A100190801, GA\v CR 201/06/0018.}
\begin{abstract}
We show that if $X$ is an infinite-dimensional separable Banach space (or more generally a Banach space with an infinite-dimensional separable quotient)
then there is a continuous mapping $f\colon X\to X$ such that the autonomous differential equation $x'=f(x)$ has no solution at any point.
\end{abstract}
\maketitle

In order to put our results into context, let us start by formulating the classical theorem of Peano.
\begin{theorem}[Peano]\label{t:peano}
Let $X=\R^n$, $f\colon\R\times X\to X$ be a continuous mapping, $t_0\in\R$, $x_0\in X$.
Then the ordinary differential equation
\begin{equation}\label{e:ode}
x'=f(t,x)
\end{equation}
together with an initial condition
\begin{equation}\label{e:inic}
x(t_0)=x_0
\end{equation}
has a solution on some open interval containing $t_0$.
\end{theorem}

Using an infinite dimensional Banach space $X=c_0$, Dieudonn\'e \cite{D} constructed a counterexample to Theorem~\ref{t:peano}.
Many counterexamples in various infinite dimensional Banach spaces followed,
e.g. \cite{LL}, \cite{B}, \cite{H}, \cite{Y}, \cite{G1}, and \cite{C} for every non-reflexive Banach space.
Finally, Godunov in \cite{G3} proved that Theorem~\ref{t:peano} is false in every infinite-dimensional Banach space.
More precisely, for every infinite-dimensional Banach space $X$, $t_0\in\R$, $x_0\in X$, there exists a continuous mapping $f\colon\R\times X\to X$,
such that there exists no solution to the initial value problem \eqref{e:ode}, \eqref{e:inic}.
Nevertheless, the above constructions are in fact showing the failure of the condition \eqref{e:inic},
and the constructed examples have many solutions on intervals not containing the given time $t_0$.

\smallskip
Moreover, Lasota and Yorke \cite{LY} (see also Vidossich \cite{V}) proved that for every Banach space $X$ and every initial
condition $t_0\in\R$, $x_0\in X$, the set of all continuous mappings $f\colon\R\times X\to X$
such that the initial value problem \eqref{e:ode}, \eqref{e:inic} has a solution that exists on the whole real line, is a generic set.
More precisely, putting the topology of uniform convergence on the space of all continuous mappings $f\colon\R\times X\to X$,
the set of mappings admitting a solution has a complement of first Baire category.

\smallskip
In view of these results, in infinite-dimensional Banach spaces it is natural to consider a weaker form of Peano's theorem:
\begin{theorem}[weak form of Peano's theorem]\label{t:peano-w}
Let $X=\R^n$, $f\colon\R\times X\to X$ be a continuous mapping.
Then the ordinary differential equation
\[x'=f(t,x)\]
has a solution on some open interval.
\end{theorem}

Showing the failure of this theorem in infinite-dimensional Banach spaces is clearly a harder problem.
In \cite{G2}, Godunov constructed a counterexample to Theorem~\ref{t:peano-w} in the Hilbert space.
Finally, Shkarin \cite{S} proved that Theorem~\ref{t:peano-w} fails for every Banach space $X$ that has a complemented subspace with an unconditional Schauder basis.
(To be more precise, Shkarin's result is even stronger
as it contains a precise quantitative information on the modulus of continuity of $f$ --- we refer the reader to \cite{S}).
Shkarin's result applies to many ``classical'' Banach spaces, such as $L_p$, $1\le p<\infty$, or $C[0,1]$.
However, there exist separable reflexive Banach spaces that contain no unconditional basic sequence (\cite{F}).
Similarly, the classical non-separable Banach space $\ell_\infty$ is also not covered by \cite{S}
(because all of its complemented subspaces are again isomorphic to $\ell_\infty$, a result of Rosenthal, see \cite{LT}).

\smallskip
The main result of this note, Theorem~\ref{cor-fin}, states that if $X$ is a Banach space with an infinite-dimensional separable quotient
(in particular every infinite-dimensional separable Banach space, of course) then Theorem~\ref{t:peano-w} fails to be true for some continuous mapping $f$.
A slightly stronger result holds, namely there is a continuous mapping $f\colon X\to X$ such that the autonomous differential equation $x'=f(x)$ has no solution at any point.
We note that the question whether every Banach space has a separable quotient is one of the outstanding problems of the Banach space theory.
It is known to hold in all reasonable classes of Banach spaces, such as reflexive, weakly compactly generated,
$C(K)$ where $K$ is a compact space (\cite{HMVZ}, Corollary~5.43, Exercise~5.8), etc.
We refer to \cite{FHHMPZ}, \cite{LT}, \cite{DGZ} and \cite{DMNZ} for background in Banach space theory and differential equations.

\bigskip
Let $\mc S$ denote a class of Banach spaces such that
$X\in\mc S$ whenever there is a continuous mapping $f\colon\R\times X\to X$ for which the equation $x'=f(t,x)$ has no solutions.
We start with two lemmata. The first one allows us to construct autonomous equations from the general ones.
\begin{lemma}\label{l:aut}
Let $X$ be a Banach space such that it has a proper complemented subspace of class $\mc S$.
Then there is a continuous mapping $f\colon X\to X$ such that the autonomous equation $x'=f(x)$ has no solutions.
\end{lemma}
\begin{proof}
The idea of the proof comes from \cite{S}.

Let $Y$ be a proper subspace of $X$ complemented by a projection $P$ such that $Y\in\mc S$.
Let $g\colon\R\times Y\to Y$ be a continuous mapping such that the equation $y'=g(t,y)$ has no solutions.
Pick $e\in\ker P$, $e\neq 0$, and find $\vp\in Y^\perp$ such that $\vp(e)=1$.
Define a mapping $f\colon X\to X$ by
\[f(x)=e+g(\vp(x),Px).\]
Obviously $f$ is continuous.

Suppose the equation $x'=f(x)$ has a solution,
i.e. there is an open interval $I$ and a mapping $x\colon I\to X$ satisfying $x'(t)=f(x(t))$ for every $t\in I$.
Then
\[\bigl(\vp(x(t))\bigr)'=\vp(x'(t))=\vp\bigl(f(x(t))\bigr)=\vp(e)=1 \quad\text{for $t\in I$,}\]
and hence $\vp(x(t))=t+c$ for $t\in I$ and some $c\in\R$.
Now define a mapping $y\colon I+c\to Y$ by $y(t)=Px(t-c)$.
Then
\[y'(t)=Px'(t-c)=Pf(x(t-c))=g\bigl(\vp(x(t-c)),Px(t-c)\bigr)=g(t,y(t)) \quad\text{for $t\in I+c$,}\]
which is a contradiction.

\end{proof}

The second lemma allows us to prove our theorem for spaces with a Schauder basis and then extend the result to much larger class.
\begin{lemma}\label{l:quot}
Let $X$ be a Banach space with a quotient from $\mc S$. Then $X\in\mc S$.
\end{lemma}
\begin{proof}
Let $Q\colon X\to Y$ be a quotient mapping such that $Y\in\mc S$ and
let $g\colon\R\times Y\to Y$ be a continuous mapping such that the equation $y'=g(t,y)$ has no solutions.

Let $\Psi\colon Y\to X$ be a continuous Bartle-Graves selector, i.e. $Q\Psi(y)=y$ for every $y\in Y$ (see \cite[Lemma~VII.3.2]{DGZ}).
Define the mapping $f\colon\R\times X\to X$ by $f(t,x)=\Psi(g(t,Qx))$. Obviously $f$ is continuous.
Moreover, if $x\colon I\to X$ is a solution of the equation $x'=f(t,x)$, then the mapping $y=Q\circ x$ satisfies
\[y'(t)=Qx'(t)=Qf(t,x(t))=Q\Psi\bigl(g(t,Qx(t))\bigr)=g(t,y(t)) \quad\text{for $t\in I$,}\]
which is a contradiction.

\end{proof}

From now on we will be dealing with Banach spaces that have a Schauder basis.
Let $P_k$, $k\in\N$, denote the canonical projections associated with a Schauder basis.
We put $P_0=0$, $R_0=I$, and $R_k x=x-P_k x$, $Q_k x=x-P_{k-1}x$ for $k\in\N$.

The next lemma is perhaps of independent interest.
It is an analogue of the estimate for the norm of perturbed vector in a space with an unconditional Schauder basis.
Notice however, that without unconditionality we are allowed only to take monotone perturbations.
\begin{lemma}\label{l:monotoneperturb}
Let $\{e_n;f_n\}$ be a Schauder basis of a Banach space $X$ and $\{\al_n\}\subset[0,1]$ be a real sequence.
Suppose that one of these conditions hold:
\begin{enumerate}[(a)]
\item $\norm{R_n}=1$ for each $n\in\N$ and $\{\al_n\}$ is non-decreasing,
\item $\{e_n\}$ is monotone and $\{\al_n\}$ is non-increasing.
\end{enumerate}
Then
\[\norm{\sum_{n=1}^\infty\al_n f_n(x)e_n}\le\norm x\quad\text{for each $x\in X$.}\]
\end{lemma}
\begin{proof}
Suppose (a) holds.
We first prove the statement under additional assumption that $\{\al_n\}$ is eventually constant.
By passing to $R_m x$ for a suitable $m\in\N$ if needed we may without loss of generality assume that $\al_1>0$.
Let $N\in\N$ be such that $\al_n=\al_N$ for all $n\ge N$.
Define $y_N=\al_N x$ and
\[y_n=R_ny_{n+1} + \frac{\al_n}{\al_{n+1}}P_n y_{n+1} \quad \text{for $N>n\ge 1$.}\]
It is easy to check by induction that $y_n=\al_n P_{n-1}x + \sum_{k=n}^N\al_k f_k(x)e_k + \al_N R_N x$ for $N\ge n\ge 1$.
Therefore $y_1=\sum_{n=1}^\infty\al_n f_n(x)e_n$.

Notice further that $y_n=\left(1-\frac{\al_n}{\al_{n+1}}\right)R_n y_{n+1} + \frac{\al_n}{\al_{n+1}}y_{n+1}$ for $1\le n< N$.
Thus, by the convexity of the norm and the fact that $0\le \frac{\al_n}{\al_{n+1}}\le 1$,
we obtain $\norm{y_n}\le\norm{y_{n+1}}$.
Since obviously $\norm{y_N}\le\norm x$, it follows that $\norm{y_1}\le\norm x$.

\smallskip
Now let $\{\al_n\}$ be an arbitrary non-decreasing sequence.
Let $z_n=\sum_{k=1}^n\al_k f_k(x)e_k + \al_n R_n x$ for $n\in\N$.
For any $m,n\in\N$, $n<m$, we have
\[z_m-z_n=\sum_{k=n+1}^m (\al_k-\al_n)f_k(x)e_k + (\al_m-\al_n)R_m x.\]
Applying the statement proven so far to the vector $R_n x$ and the sequence $\beta_k=\al_k-\al_n$ for $n<k\le m$,
$\beta_k=0$ for $k\le n$, and $\beta_k=\al_m-\al_n$ for $k>m$, we obtain $\norm{z_m-z_n}\le\norm{R_n x}$.
It follows that $\{z_n\}$ is convergent.
Clearly, $\lim z_n=\sum_{n=1}^\infty\al_n f_n(x)e_n$, and as $\norm{z_n}\le\norm x$ by the first part of the proof, the statement follows.

The proof under the assumption (b) is similar.

\end{proof}

Next we will define some mappings useful for our construction.
For $r>0$ we define a function $\vp_r\colon\R\to\R$ by
\[\vp_r(t)=\begin{cases}
0 & \text{for $t\in(-\infty,-2r]$,}\\
2+t/r & \text{for $t\in[-2r,-r]$,}\\
1 & \text{for $t\in[-r,r]$,}\\
2-t/r & \text{for $t\in[r,2r]$,}\\
0 & \text{for $t\in[2r,+\infty)$.}
\end{cases}\]

Suppose $X$ is a Banach space with a Schauder basis $\{e_k;f_k\}$.
We define a mapping $\Phi_r\colon X\to X$ by
\[\Phi_r(x)=\sum_{k=1}^\infty \vp_r\bigl(\norm{Q_k x}\bigr)f_k(x)e_k.\]
Notice that the mapping is well-defined as $\vp_r\bigl(\norm{Q_k x}\bigr)=1$ for all $k$ large enough.
This mapping was already used by Shkarin in \cite{S};
however, lacking the unconditionality we have to employ Lemma~\ref{l:monotoneperturb} to prove some properties of this mapping.

\begin{lemma}\label{l:Phiprop}
Let $X$ be a Banach space with a Schauder basis $\{e_k;f_k\}$ satisfying $\norm{R_n}=1$ for each $n\in\N$ and let $r>0$.
Then the mapping $\Phi_r$ has the following properties:
\begin{enumerate}[(i)]
\item\label{PhipropQcomm} $Q_n\Phi_r(x)=\Phi_r(Q_n x)$ for each $x\in X$ and $n\in\N$.
\item\label{PhipropTailid} If $\norm{Q_N x}\le r$ for some $x\in X$ and $N\in\N$, then $Q_n\Phi_r(x)=Q_n x$ for all $n\ge N$.
\item\label{PhipropLenorm} $\norm{\Phi_r(x)}\le\norm x$ for each $x\in X$.
\item\label{PhipropBdd} $\norm{\Phi_r(x)}<2r$ for each $x\in X$.
\item\label{PhipropCont} $\Phi_r$ is continuous.
\end{enumerate}
\end{lemma}
\begin{proof}
(\ref{PhipropQcomm}):
Choose $x\in X$ and compute
\[Q_n\Phi_r x=\sum_{k=n}^\infty \vp_r\bigl(\norm{Q_k x}\bigr)f_k(x)e_k=\sum_{k=1}^\infty \vp_r\bigl(\norm{Q_kQ_n x}\bigr)f_k(Q_n x)e_k=\Phi_r(Q_n x).\]

(\ref{PhipropTailid}):
Suppose $\norm{Q_N x}\le r$ and $n\ge N$. Then $\norm{Q_k x}\le r$ for all $k\ge N$ and hence $\vp_r\bigl(\norm{Q_k x}\bigr)=1$ for all $k\ge N$.
Thus
\[Q_n\Phi_r(x)=\sum_{k=n}^\infty 1\cdot f_k(x)e_k=Q_n x.\]

(\ref{PhipropLenorm}):
Since the sequence $\left\{\vp_r\bigl(\norm{Q_k x}\bigr)\right\}\subset[0,1]$ is non-decreasing we may apply Lemma~\ref{l:monotoneperturb}.

(\ref{PhipropBdd}):
Pick any $x\in X$. Let $n\in\N$ be the smallest such that $\norm{Q_n x}<2r$.
Then, using (\ref{PhipropQcomm}) and (\ref{PhipropLenorm}),
\[\norm{\Phi_r(x)}=\norm{\sum_{k=n}^\infty \vp_r\bigl(\norm{Q_k x}\bigr)f_k(x)e_k}=\norm{Q_n\Phi_r(x)}=\norm{\Phi_r(Q_n x)}\le\norm{Q_n x}<2r.\]

(\ref{PhipropCont}):
Fix any $x\in X$ and $\ve>0$. There is $n\in\N$ such that $\norm{R_n x}<\frac\ve4$.
By the continuity of the projection $R_n$ there is a neighbourhood $U$ of $x$ such that $\norm{R_n y}<\frac\ve4$ for any $y\in U$.
Further, choose a neighbourhood $V$ of $x$, $V\subset U$, such that
$\abs{\vp_r\bigl(\norm{Q_k x}\bigr)f_k(x)-\vp_r\bigl(\norm{Q_k y}\bigr)f_k(y)}\norm{e_k}<\frac\ve{2n}$ for any $y\in V$ and all $1\le k\le n$.
This can be done using the continuity of the mappings involved.
Then, using (\ref{PhipropQcomm}) and (\ref{PhipropLenorm}), we obtain for any $y\in V$ that
\[\begin{split}
\norm{\Phi_r(x)-\Phi_r(y)}&\le\norm{P_n\Phi_r(x)-P_n\Phi_r(y)}+\norm{R_n\Phi_r(x)}+\norm{R_n\Phi_r(y)}\\
&\le\sum_{k=1}^n\absb{Big}{\vp_r\bigl(\norm{Q_k x}\bigr)f_k(x)-\vp_r\bigl(\norm{Q_k y}\bigr)f_k(y)}\norm{e_k}+\norm{\Phi_r(R_n x)}+\norm{\Phi_r(R_n y)}\\
&<\frac\ve2+\norm{R_n x}+\norm{R_n y}<\ve.
\end{split}\]

\end{proof}

\begin{theorem}\label{t:basis}
Let $X$ be an infinite-dimensional Banach space with a Schauder basis $\{e_k;f_k\}$.
Then $X\in\mc S$.
\end{theorem}
\begin{proof}
Without loss of generality we may and do assume that $\norm{R_k}=1$ for each $k\in\N$.

Choose an $a\in X$ such that $Q_k a\neq 0$ for any $k\in\N$.
Define a mapping $h\colon\R\times X\to X$ by
\[h(t,x)=\begin{cases}
2t\Phi_1\left(\frac x{t^2}\right) & \text{for $t>0$,}\\
0 & \text{for $t=0$,}\\
2t\Phi_1\left(\frac{x-a}{t^2}\right) & \text{for $t<0$.}
\end{cases}\]
The mapping $h$ is continuous.
Indeed, the continuity at points $(t,x)$ for $t\neq 0$ follows from the continuity of the mapping $\Phi_1$,
while the continuity at points $(0,x)$ follows from the boundedness of the mapping $\Phi_1$.
We note that the mapping $h$ has the property that the equation $x'=h(t,x)$ has no solutions on any interval containing zero. (See also \cite{S}.)
Notice also that by Lemma~\ref{l:Phiprop}(\ref{PhipropBdd})
\begin{equation}\label{e:hestimate}
\norm{h(t,x)}\le 4\abs t \quad\text{for all $t\in\R$, $x\in X$.}
\end{equation}

In \cite{S} Shkarin constructs an equation with no solutions by splitting the space into countably many infinite-dimensional pieces
and using a copy of $h$ in each of these pieces shifted to a different time.
This splitting is impossible without using the unconditionality and hence we have to develop a different approach.
We will spread copies of the mapping $h$ directly over the time axis.

Choose a sequence $\{\ve_n\}$ such that $\ve_n>0$ for each $n\in\N$ and
\begin{equation}\label{e:sumeps}
\sum_{n=1}^\infty\ve_n<1.
\end{equation}
Let $\{t_n\}$ be an enumeration of rational numbers such that $t_i\neq t_j$ for $i\neq j$.

By induction, for each $n\in\N$ we find numbers $\de_n$, $u^n_i$, $v^n_i$ for $i\in\N\cup\{-1,0\}$, satisfying the following conditions:
\begin{enumerate}[(D1)]
\item\label{cndIrat} $\de_n>0$ and $\de_n\in\R\setminus\Q$.

\item\label{cnd1} $8\de_n\le\ve_n/2$.

\item\label{cndInside} For each $k\in\{1,\dotsc,n-1\}$ the following holds:
If $t_n\in(u^k_j,u^k_{j+1})$ for some $j\in\N\cup\{-1,0\}$ then $2\de_n<t_n-u^k_j$ and $2\de_n<u^k_{j+1}-t_n$.
If $t_n\in(v^k_{j+1},v^k_j)$ for some $j\in\N\cup\{-1,0\}$ then $2\de_n<v^k_j-t_n$ and $2\de_n<t_n-v^k_{j+1}$.

\item\label{cndMeas} For each $k\in\{1,\dotsc,n-1\}$ the following holds:
Suppose $t_n\in(u^k_j,u^k_{j+1})\cup(v^k_{j+1},v^k_j)$ for some $j\in\N\cup\{0\}$. Then
\[4\de_n<\frac1{2^n}\frac1{2^{j+2}}\left(\left(\frac45\right)^{j+1}\de_k\right)^2.\]

\item $u^n_i=t_n-\left(\frac45\right)^i\de_n$ and $v^n_i=t_n+\left(\frac45\right)^i\de_n$ for $i\in\N\cup\{0\}$, $u^n_{-1}=-\infty$, $v^n_{-1}=+\infty$.
\end{enumerate}

Finally we define $g_0(t,x)=0$ and
\[g_n(t,x)=g_{n-1}(t,x)+\vp_{\de_n}(t-t_n)\left(h(t-t_n,x)-\Phi_{\ve_n/4}(g_{n-1}(t,x))\right)\]
for all $t\in\R$, $x\in X$, $n\in\N$

All the mappings $g_n$ are continuous since both $h$ and $\Phi_{\ve_n/4}$ are continuous.
Further, using \eqref{e:hestimate}, Lemma~\ref{l:Phiprop}(\ref{PhipropBdd}), the properties of function $\vp_{\de_n}$, and condition (D\ref{cnd1}) we obtain for any $n\in\N$ that
\[\norm{g_n(t,x)-g_{n-1}(t,x)}<8\de_n+\frac{\ve_n}2\le\ve_n \quad\text{for all $t\in\R$, $x\in X$.}\]
It follows that the sequence of mappings $\{g_n\}$ converges uniformly on $\R\times X$ to a continuous mapping $g\colon\R\times X\to X$.
We claim that the equation $x'=g(t,x)$ has no solutions.

Let us prove it by contradiction.
Suppose there is an open interval $I\subset\R$ and a mapping $y\colon I\to X$ such that $y'(t)=g(t,y(t))$ for each $t\in I$.
Find $N\in\N$ such that $t_N\in I$. To simplify our notation we denote $u_i=u^N_i$, $v_i=v^N_i$ for $i\in\N$.

From the properties of the functions $\vp_{\de_k}$ it follows that
\begin{equation}\label{e:geqgn}
g(t,x)=g_N(t,x) \quad\text{for each $x\in X$ and $t\in\R\setminus\Delta$,}
\end{equation}
where $\Delta=\bigcup\limits_{k=N+1}^\infty (t_k-2\de_k,t_k+2\de_k)$.
Denote $\Delta_j=\Delta\cap(u_j,u_{j+1})$ for $j\in\N$.
Notice that by (D\ref{cndIrat}) and (D\ref{cndInside}) we have $\Delta_j=\!\!\!\bigcup\limits_{\substack{t_k\in(u_j,u_{j+1}) \\ k>N}}\!\!\!(t_k-2\de_k,t_k+2\de_k)$.
Thus by (D\ref{cndMeas})
\begin{equation}\label{e:Djmeasureest}
\lambda(\Delta_j)\le\!\!\sum_{\substack{t_k\in(u_j,u_{j+1}) \\ k>N}}\!\!\!\!4\de_k<\sum_{k=N+1}^\infty \frac1{2^k}\frac1{2^{j+2}}\left(\left(\frac45\right)^{j+1}\de_N\right)^2<\frac1{2^{j+2}}\left(\left(\frac45\right)^{j+1}\de_N\right)^2,
\end{equation}
where $\lambda$ denotes the Lebesgue measure on $\R$.

Further, for each $x\in X$ and $t\in\R$ such that $\abs{t-t_N}<\de_N$ there is $U_{t,x}$ a neighbourhood of $(t,x)\in\R\times X$ and $K_{t,x}\in\N$ such that
\begin{equation}\label{e:gnloch}
Q_k g_N(s,y)=Q_k h(s-t_N,y) \quad\text{for any $k\ge K_{t,x}$ and $(s,y)\in U_{t,x}$.}
\end{equation}
Indeed, find $K_{t,x}\in\N$ so that $\normb{big}{Q_{K_{t,x}}g_{N-1}(t,x)}<\ve_N/4$
and utilising the continuity of $g_{N-1}$ choose $U_{t,x}$ in such a way that $(s,y)\in U_{t,x}$ implies $\normb{big}{Q_{K_{t,x}}g_{N-1}(s,y)}<\ve_N/4$
and moreover $\abs{s-t_N}<\de_N$. Thus for $(s,y)\in U_{t,x}$ we have $\vp_{\de_N}(s-t_N)=1$ and hence for $k\ge K_{t,x}$ by Lemma~\ref{l:Phiprop}(\ref{PhipropTailid})
\[Q_k g_N(s,y)=Q_k g_{N-1}(s,y)+Q_k h(s-t_N,y)-Q_k \Phi_{\ve_N/4}(g_{N-1}(s,y))=Q_k h(s-t_N,y).\]

Let $i\in\N$ be such that $u_i\in I$ and $v_i\in I$.
Since the set $W=\{(t,y(t));\;t\in[u_i,v_i]\}\subset\R\times X$ is compact, there is a finite subcovering $\{U_j\}_{j=1}^m$ of a covering $\{U_{t,y(t)};\;t\in[u_i,v_i]\}$ of $W$.
Denote the constants from \eqref{e:gnloch} corresponding to $U_j$ by $K_j$, $j=1,\dotsc,m$.
Find $K\in\N$ such that $K\ge\max_{1\le j\le m}K_j$,
\begin{equation}\label{e:yleonehalf}
\norm{Q_K(y(u_i)-a)}<\frac12(t_N-u_i)^2, \quad\text{and}\quad \norm{Q_K y(v_i)}<\frac12(v_i-t_N)^2.
\end{equation}
From \eqref{e:gnloch} it follows that
\begin{equation}\label{e:geqh}
Q_K g_N(t,y(t))=Q_K h(t-t_N,y(t)) \quad\text{for every $t\in[u_i,v_i]$.}
\end{equation}

We claim that
\begin{equation}\label{e:yincone}
\norm{Q_K(y(t)-a)}<(t_N-t)^2 \quad\text{ for every $t\in[u_i,t_N)$.}
\end{equation}
We prove this by induction, showing that for $j\ge i$
\begin{equation}\label{e:yleonuj}
\norm{Q_K(y(t)-a)}<\left(1-\frac1{2^{j+1}}\right)(t_N-t)^2 \quad\text{for every $t\in[u_j,u_{j+1}]$.}
\end{equation}
So let $j\ge i$ and suppose $\norm{Q_K(y(u_j)-a)}<(1-\frac1{2^j})(t_N-u_j)^2$ from the previous induction step (or from \eqref{e:yleonehalf} for the first step).
Assume \eqref{e:yleonuj} does not hold and put
\begin{equation}\label{e:yleavescone}
S=\inf\left\{t\in[u_j,u_{j+1}];\;\norm{Q_K(y(t)-a)}\ge\left(1-\frac1{2^{j+1}}\right)(t_N-t)^2\right\}.
\end{equation}
Then $S\in(u_j,u_{j+1}]$.
Define $z(t)=Q_K a+\frac{Q_K(y(u_j)-a)}{(t_N-u_j)^2}(t_N-t)^2$ for $t\in\R$ and find $T\in[u_j,S]$ for which
\begin{equation}\label{e:Tmaxdist}
\norm{Q_K y(T)-z(T)}=\max_{t\in[u_j,S]}\norm{Q_K y(t)-z(t)}.
\end{equation}
We have $z'(t)=\frac{2(z(t)-Q_K a)}{t-t_N}$ for $t\neq t_N$ and hence
\[z(T)=z(u_j)+\int_{u_j}^T \frac{2(z(t)-Q_K a)}{t-t_N}\d t=Q_K y(u_j)+\int_{u_j}^T \frac{2(z(t)-Q_K a)}{t-t_N}\d t.\]
Since $y'(t)=g(t,y(t))$ for each $t\in I$, we have
\[y(T)=y(u_j)+\int_{u_j}^T y'(t)\d t=y(u_j)+\int_{u_j}^T g(t,y(t))\d t\]
and hence
\[Q_K y(T)=Q_K y(u_j)+\int_{u_j}^T Q_K g(t,y(t))\d t.\]
To prove our claim we show that $Q_K y$ does not deviate too much from $z$, which is a solution to a ``non-perturbed'' equation.

Using \eqref{e:geqh} and then Lemma~\ref{l:Phiprop}(\ref{PhipropTailid}) together with \eqref{e:yleavescone} we obtain
$Q_K g_N(t,y(t))=2(t-t_N)Q_K\Phi_1\left(\frac{y(t)-a}{(t-t_N)^2}\right)=\frac{2Q_K(y(t)-a)}{t-t_N}$ for each $t\in[u_j,S]$.
Thus
\[\begin{split}
Q_K y(T)-z(T)&=\int_{u_j}^T\left(Q_K g(t,y(t))-\frac{2(z(t)-Q_K a)}{t-t_N}\right)\d t\\
&=\int_{u_j}^T Q_K\bigl(g(t,y(t))-g_N(t,y(t))\bigr)\d t+\int_{u_j}^T\left(Q_K g_N(t,y(t))-\frac{2(z(t)-Q_K a)}{t-t_N}\right)\d t\\
&=\!\!\!\int_{(u_j,T)\cap\Delta}\!\!\!\!\! Q_K\bigl(g(t,y(t))-g_N(t,y(t))\bigr)\d t+2\int_{u_j}^T \frac1{t-t_N}\bigl(Q_K y(t)-z(t)\bigr)\d t,
\end{split}\]
where we use also \eqref{e:geqgn} to verify the last equality.
Applying \eqref{e:sumeps} and \eqref{e:Tmaxdist} we can estimate
\[\begin{split}
\norms{Q_K y(T)-z(T)}&\le\int_{\Delta_j} 1\d t+2\norms{Q_K y(T)-z(T)}\int_{u_j}^{u_{j+1}} \frac1{t_N-t}\d t=\lambda(\Delta_j)+2\norms{Q_K y(T)-z(T)}\log\frac{t_N-u_j}{t_N-u_{j+1}}\\
&=\lambda(\Delta_j)+2\norms{Q_K y(T)-z(T)}\log\frac54<\lambda(\Delta_j)+\frac12\norms{Q_K y(T)-z(T)}
\end{split}\]
and so finally by \eqref{e:Tmaxdist} we obtain
\[\norm{Q_K y(S)-z(S)}\le\norms{Q_K y(T)-z(T)}<2\lambda(\Delta_j).\]
But this last inequality together with \eqref{e:Djmeasureest} implies
\[\begin{split}
\norm{Q_K(y(S)-a)}&\le\norm{Q_K y(S)-z(S)}+\norm{z(S)-Q_K a}<2\lambda(\Delta_j)+\frac{\norm{Q_K(y(u_j)-a)}}{(t_N-u_j)^2}(t_N-S)^2\\
&<2\lambda(\Delta_j)+\left(1-\frac1{2^j}\right)(t_N-S)^2<\frac1{2^{j+1}}\left(\left(\frac45\right)^{j+1}\de_N\right)^2+\left(1-\frac1{2^j}\right)(t_N-S)^2\\
&=\frac1{2^{j+1}}(t_N-u_{j+1})^2+\left(1-\frac1{2^j}\right)(t_N-S)^2\le\left(1-\frac1{2^{j+1}}\right)(t_N-S)^2,
\end{split}\]
a contradiction with~\eqref{e:yleavescone}.
This finishes the proof of~\eqref{e:yincone}.

Now \eqref{e:yincone} immediately implies that $\lim\limits_{t\to t_N-}Q_K y(t)=Q_K a\neq 0$.
Further, analogously as above (approaching the point $t_N$ from the right, replacing $u_j$s by $v_j$s and $a$ by $0$) we can show that $\lim\limits_{t\to t_N+}Q_K y(t)=0$.
These two facts contradict the continuity of $y$ at $t_N$.

\end{proof}

As every separable infinite-dimensional Banach space has an infinite-dimensional quotient space with a Schauder basis (see \cite[Theorem~1.b.7]{LT}),
Theorem~\ref{t:basis} together with Lemma~\ref{l:quot} and Lemma~\ref{l:aut} give us the final result:
\begin{theorem}\label{cor-fin}
Let $X$ be a Banach space with an infinite-dimensional separable quotient.
Then there is a continuous mapping $f\colon X\to X$ such that the autonomous equation $x'=f(x)$ has no solutions.
\end{theorem}

\end{document}